\documentclass[11pt,twoside]{article}
\usepackage{amsfonts}
\begin{document}
\centerline{\large{\bf Embodied Cognition and the Origins of
Geometry:}}
 \centerline{\large{\bf A Model Approach of Embodied Mathematics}}
\centerline{\large{\bf Through Geometric Considerations }}

\par \bigskip
{\bf Dionyssios Lappas},
\par
Department of Mathematics,
\par
University of Athens,
\par
Panepistimiopolis 157 84 Athens, Greece
\par
tel. 210 7276394 FAX: 210 727 6378   dlappas@cc.uoa.gr
\par \bigskip
 {\bf Panayotis Spyrou},
 \par
Department of Mathematics,
\par
University of Athens,
\par
Panepistimiopolis 157 84 Athens, Greece
\par
tel. 210 7276392  FAX: 210 727 6410  pspirou@cc.uoa.gr
\par
\bigskip
\begin{abstract}
In this paper, we propose that 'embodied mathematics' should be
studied not only by reduction to the present individual bodily
experience but in an historical context as well, as far as the
origins of mathematics are concerned. Some early mathematical
results are the Theorems of Geometry and arose as attempts to
objectively render the main perceptual categories such as
verticality, horizontality, similarity (or its varieties).
Inasmuch as these are of a qualitative nature, it was required
that they be expressed in a quantitative way in order to be
objectified.
 The first form of this objectification occurred in the case of 'archetypal results', namely the Pythagorean triads and the internal ratio of the legs in the right triangles. In the next stage, a 'scientific' treatment would come from a shift of objectification and descriptions inside an abstract theory, which would constitute the first logicomathematical knowledge. In this theory, the 'archetypal results' were incorporated, generalized and acquired their unquestionable, supertemporal validity.
The study presents a particular epistemological analysis of some
of the main terms used in the beginnings of Geometrical Thought
and Euclid's Elements, utilizing the theoretical apparatus of the
theory of 'embodied mathematics'. It also traces models of
objectification for the 'archetypal results' and indicates their
diffusion in later mathematical developments.
\par
\bigskip
Key words: Archetypal results, conceptual categories, embodied
mathematics, objectification, pattern, perceptual categories,
prototypes.
\end{abstract}
\bigskip
\par

\centerline{\bf 1. The embodied situation}
\bigskip
The idea of embodied mathematics was suggested recently by Nunez
\& al. (1999) and Lakoff \& Nunez (2000). The relevance of the
human body to its conceptual systems is decisive in the theory of
embodied mind in general, as a theory in which the correlation
between human experience and cognitive sciences is attempted,
Varelas \& al. (1992). The leading idea in embodied mathematics is
the mediation of spatiality in the formation of logical thought.
In this respect, the comprehension of space (or its models) by
humans is of primary importance in the history of knowledge and
in cognition as well. Lakoff discusses embodied mathematics at
least since 1987, when he stated "mathematics is based on
structures within the human conceptual system, structures that
people used to comprehend ordinary experience" (Lakoff 1987, p.
364) and he  expressed the ambition of unraveling the mystery of
Platonic ideas (ibid.). Today, modifications of the same idea are
also used by researchers in education, Gray \& Tall (2000),
Edwards (1998), Boero P. \& Bazzini L. \& Garutti R, (2002),
Watson (2002), Watson \& Spyrou \& Tall (2003), etc. The critique
of the above idea by Presmeg (2002) or Schiralli \& Sinclair
(2003) is of great interest as well. Lakoff \& Nunez, in
particular, suggest a list of properties of embodied mathematics,
supporting the idea that mathematics is a product of the full
range of human experience:

Human mathematics is embodied, it is
grounded in bodily experience in the world... not purely
subjective... not a matter of mere social agreement...

 It uses
the very limited and constrained resources of human biology and
is shaped by the nature of our brains, our bodies, our conceptual
systems, and the concerns of human societies and cultures, Lakoff
\& Nunez  ( pp 348 -365).

On the other hand, embodied mathematics allows us to address the
main epistemological problem of mathematics (Lakoff 1987, pp. 353-369).
 We prefer the formulation of Piaget[1]:

The epistemology of mathematics has three principal and classic
problems: why is it so fruitful though based on very few and
relatively poor concepts or axioms; why has it a necessary
character, thus remaining constantly rigorous despite its
constructive character which could be a source of irrationality;
and why does it agree with experience or physical reality in
spite its completely deductive nature?, Piaget (1972, p 69).

The origin of such ideas in fact goes back to the
phenomenological tradition of Derrida (1962), Merleau-Ponty
(1954) or even that of Husserl (1917). A part of this heritage is
acceptable by Varelas \& al (1992) or Nunez \& al (1999) as well.
We do not deal with the various interpretations of the notion or
the evolution of the embodied cognition inasmuch as it is very
well described in the above respective literature. In the case of
geometry, we suggest a model of exposition for this line of ideas
by using historical references. This point of view seems to be
very little discussed and developed so far.

 In this regard and in
order to describe the embodied situation, we have to reflect on
perception, not only in the context of contemporary neuroscience
but also in the light of historical epistemology. This genetic
approach inquires into the potential environments, the necessary
metaphors and ideas that have preceded the scientific formation
in mathematics, as we know it today.

It is remarkable that, in
Husserl's writings we come across a program for a recursive
inquiry into the "origins of geometry", as the basis of our
culture [2]. This program intended to discover the ground of the
axioms in the pre-scientific period of geometry, since the
established symbolism has covered over the origins. However, the
human capacity of using and manipulating symbols has become
dominant in the evolution of civilization. But, if we deal with
geometry, special attention should be given to the fact that the
very symbols used in this domain (such as point, line, etc.) must
resemble the mental object we seek to symbolize (Beth \& Piaget,
p. 217). In this sense, we view these correspondences (between
geometric symbols and concepts) as inextricably bound to our
bodily experience and we advocate that embodied mathematics can
be presented in an historical perspective, especially as far as
their origins, which are related to geometry, are concerned.
\par
\bigskip
\centerline{\bf 1.1 Perception and the formation of geometric
notions}
\bigskip
In order to investigate the embodied primitive concepts, another
web of ideas is needed that connects us with our hypostasis as
intelligent beings and is inherent in our perceptual bodily
experience in a nature dominated by gravity [3]. Note that
Merleu-Ponty underlined the significance of gravity in the
comprehension of space and his conviction is supported by others,
ever since: Zuzne, (1970), Ibbotson \& Bryant (1976), Varelas \&
al (1992) and Lakoff (1987), indicatively. Actually, we need to
distinguish the main components of this perceptual system insofar
as they seem to be related to geometry. According to Lakoff \&
Johnson (1980 p. 277) verticality is the main source domain and
is connected with our understanding of quantity as well. There is
also evidence that the function of apprehension of the vertical
has an affinity to the horizontal, as has been proposed in
psychology; the acquisition of the vertical is synchronous with
that of the horizontal, Piaget \& Inhelder (1956), Mackay \& al.
(1972). Apart from these two perceptual categories, another
important aspect of the cognitive apparatus is the perception of
Gestalts [4] that gives us the ability to recognize the 'shape' or
'form' of objects.

 Subsequently, we will focus on the above three
functions that mainly constitute our perceptual apparatus, i.e.,
the vertical, the horizontal, and the recognition of the shape (:
our ability to apprehend the similarity of shapes). Altogether
they are connected functionally, decisively affect our adaptation
in the environment, and offer us the main mental tools to
describe our experiences in the world. We should note here that
similarity generally is a far broader perceptual category: it
appears either as visual (figures or coloring) or auditory etc,
and consists of a general trend both the perceptual and
conceptual systems that unify the manifold of experience into
rules, Gentner \& Medina (1998). Thus the recognition of visual
shape appears to be a particular function of the skill of
apprehending similarity. If we consider similarity as a quality,
it is difficult to communicate since qualities are subjective
impressions. Finally, by using numbers to interpret aspects of
similarity and to communicate them; thereby, a sort of inter-
subjective knowledge is achieved [5].

 In this paper we claim
something stronger; human bodily experience of the world is
transformed to rationality through the mediation of the geometric
comprehension of the world [6]. In particular, we claim that the
first results of Geometry arose from the persistent effort for an
objective rendering and reification of the three main perceptual
categories of space forms, via arithmetic and even logical
relationships.  In their first formulation such results  (namely
the Pythagorean triads, the internal ratio of the legs in the
right triangles, etc) occurred as 'archetypal' formations of
later basic theorems of Geometry. On the other hand, along these
lines, a reduction of the formation and development of the
geometric thought to a psychological background is also suggested.
\par
\bigskip
\centerline{\bf 2. Archetypal results, Patterns and Prototypes in
Geometry}
\bigskip
    Mathematics is founded on the logicomathematical ability
     of enumeration. Humankind's apprehension of the idea of number was a
      cornerstone of objectivity. The significance of numbers increased when
       numerical systems improved and allowed the enumeration of huge quantities,
        and made their operations simpler. Through numbers, the experience becomes
         homogeneous, inter-subjective easier to transmit and as a result, number
          stands for a certain grasping of the world's truth [7].
 Another set of practical questions (connected to Geometry), in which
 numbers were involved, had to do with the measurement of land areas.
Beyond this task, the construction of homes, temples and pyramids evoked deep
 questions related to the classification of shapes. In order to understand
these spatial phenomena and especially their contribution to the formation of
 geometric concepts Bender and Schreiber [8] (1980) suggested the very idea
  of norm and argued in favor of linking the connection of daily activities
   and notions of geometric shapes.

At any rate, using current terms, we suggest another exposition of the
above process. In this, a series of brain and linguistic functions are
intervene in our conceptual systems, mainly of prototypes,
a linguistic term introduced by Rosch (1978). In our context a prototype is
 regarded as the 'best exemplar' of a concept (and is often
  'a non-existence'), as a specification for the members of a family and in
   this sense, it is the ideal core of a concept, Harley (1995, p. 193),
    Malt [9] (1999, p. 333). Lines, points, and planes, in their ideal use,
     serve as our main examples for prototypes in Geometry.
The prototypes are also involved in the conceptual evolution of
student's geometric thought Tall (1995), Tall \& al (2000). Tall
\& al attempted a weak description of the involvement of
prototypes in developing cognitive structures and advanced
mathematical thinking.

In our study, we suggest that between the
primary notions articulated via prototypes and the archetypal
geometric results an additional conceptual structure should
mediate. This process is carried out by means of a constructed
hierarchy referring to notions of growing complexity, for which
we propose the term patterns. In this hierarchy first comes the
pattern of a triangle, as it involves a minimum of data (point,
line and plane) in a coherent structure. (Note that next in this
classification could be the tetrahedron, a distinct spatial
pattern).

 As a pattern, the triangle is the first step for the
analysis of a figure and its reduction to prototypes. In a
triangle, the broader categories of equality, similarity, and
area, bestow a quality of additional structure and bring forth
associations (of plane geometry, for example). Thus, the role of
a triangle is fundamental and it becomes the principal
instrument, mediating in all proofs concerning more complicated
figures. According to this point of view the pattern of the
triangle is a necessary and decisive element for the later
development of Geometry.
\par
\bigskip
\centerline{\bf 3. The reduction of the archetypal results to Prototypes}
\bigskip
The archetypal results constituted potentially an objective
knowledge, before the emergence of any coherent logical deductive
theory (as is also indicated by a position of Lakatos, 1997). In
this respect, it is altogether interesting to locate the
necessary regressive processes the mind has followed in order to
refine for itself the fundamental terms and finally form a
logical theory where the archetypal results could be placed. In
any case, the meaning and the semantics of abstract geometric
thought consist of exactly these recursive functions, which lead
to the determination of the minimum and sufficient terms and
principles (the axioms) that constitute the theory [10].

 A key
instrument of investigation of Greek mathematics was focused
around the fundamental (platonic) demand of the anti-visual, as
Szabo has noticed:

It seems that new kinds of proof appeared at the same time as
Greek mathematics was becoming anti-empirical and anti-visual,
Arprad Szabo (1978, p 197).

 In order for the archetypal results to come under the control
 of logic and the requirements of the anti-visual, the ideal configuration
  that determines their possible rendering in a necessary and unique way had
  to be invented. In Euclidean Geometry for instance, we have an apparent
   result: "Two straight lines may intersect at one, and only one, point".
    Such a formulation of a property of geometric objects could only be
    arrived at after grasping the prototypical notion of straight line [11].

This approach, and the attendant type of reasoning, actually
designates a shift in the process of the objectification of
knowledge:

 But to make an object of something, to make it a subject of a predication
 or attributions, merely differs in name from having a presentation of it, and
  having a presentation in a sense which, while not the only one, is none the
   less the standard one for logic.
 (Husserl, Logical Investigation II, p.  366)
 \par
\bigskip
 \centerline{\bf   3.1 Angle: a compound Prototype}
\bigskip
The 'angle' is crucially involved in the fundamental perceptual
categories, as they are presented in Section 2, inasmuch as it
obviously mediates in shape recognition. Also the same prominent
role for the angle is reserved by mathematicians in the early
foundational processes, and thus it turns out to be an
unavoidable term in any kind of reasoning. In Euclid's Elements
(I, Definition 8 and 9), the angle is defined as follows:

A plane angle is the inclination to one another of two lines in a
plane which meet one another and do not lie in a straight line.
And when the lines containing the angle are straight, the angle
is called rectilinear. (quoted in Heath, 1956 p. 176).

Thus, at least within the scientific status of the theory, the
angle is not given as a measurable magnitude. It rather occurs as
a compound prototype, reduced to the simpler basic prototypes of
the point and straight line. The compound prototype is
anticipated by the theory [12]. The angles are only amenable to the
single manipulation of superposition and coincidence, in the same
uniform way that occurs for all figures in Euclid's Geometry.
These synthetic conditions do not constitute objective criteria
justifying the angle, independent of the subjective experience,
like those provided by measurement. Besides, two angles are
similar if and only if they are congruent. Thus the angle becomes
a form which determines the shape [13].

Note that an acceptable
mathematical measurement of the angle is not at all obvious and
turns out to be a deep result, both in mathematics [14] and in
neurology [15]. However, through a sophisticated treatment in Euclid's
Elements, the angles are subject to manipulation in an absolute
way, without any mediation of measurement or theory.

In the case
of perpendicular lines and in parallels as well, either the
definitions or their verification is reduced to certifications
concerning prototypical configurations, such as point and
straight line:

In fact Euclid does not use any relation which is not reducible
to coincidability between lines until he treats ratios. (Mueller
1981, p. 41).

In particular: A straight line meets another straight line
perpendicularly, when the formed angles are equal - then the
angles are called right angles - All right angles are equal,
(Elements I, Definition 10 and Postulate 4). In the above
statements a twofold purpose is achieved: First, it serves as the
identification condition for perpendicularity, via right angles.
Further, through the use of 'all' in the above assertion, an
implicit transcendental definition is formulated.
\par
\bigskip
\centerline{Figure 1: Lines l and l are perpendicular iff $a = b
 \perp $}
\bigskip
 In the
case of parallel lines, their definition includes a
transcendental term, such as 'produced indefinitely' (Heath,
1956, p. 190). We remark again that the formulation of the 5th
Postulate (ibid. p. 202) permitted a finitistic type of argument,
referring again to straight lines (see Figures 1 and 2).
\par
\bigskip
\centerline{Figure2: Lines l and l are parallel iff $a + b$ would
 form a straight line}
\bigskip
We consider all the above as indicative and exemplary cases of
reduction to prototypes. In fact, it is enough to ascertain that
after a suitable 'displacement', we arrive at an arrangement that
proves to be a straight angle.
\par
\bigskip
\centerline{\bf 3.2 The triangle as a Pattern}
\bigskip
As we have already noticed the triangle is introduced not only as
a figure but also as a structure, that is followed immediately by
a number of connotations, leading to the conditions of identity
and difference that will determine its ontology.   Subsequently,
a triangle - and as a consequence any complicated geometrical
object - would be acceptable if and only if we had rigidity
criteria for the particular structure that would make it
recognizable in all its potential appearances and permissible
transformations.

 This logical determination of the triangle was
set up after the establishment of the 'congruence' and
'similarity' criteria. The criteria contain some potential
actions that are advanced in mental acts of comparisons. These
universal properties provide a status for the notion of the
triangle, for which we propose the term pattern.

 The triangle
would eventually constitute the touchstone ensuring logical
comparisons and arguments. Situations of this kind arise through
the introduction of new concepts, as well as appear in the
majority of the proofs. Think, for instance, of the result
concerning the angles - sum for a triangle and its variations
that determines the type of Space Geometry (Euclidean, Spherical
or Hyperbolic).
\par
\bigskip
\centerline{\bf {4. From archetypal results to geometrical
theorems}}
\bigskip
\centerline{\bf {4.1 Pythagorean Theorem}}
\bigskip
The constitution of a theory for Geometry and in particular the
need to encompass the already existing 'archetypal' results posed
certain methodological problems. This necessitated the further
derivation of mental instruments, as the accumulated empirical
knowledge was not enough.

Concerning the Pythagorean Theorem,
such a new concept is that of a square which  denotes the second
power of a number and at the same time is aprehended as the area
of a rectangle. This situation might be considered as an early
twofold representative expression, both of arithmetic and
geometric nature, concerning two-dimensional objects. Another
cognitive factor arising is the need for formulations and proofs
in terms of the theory under development and at the same time
obeying Logic. In the case of a square, its existence seems to be
altogether not self-evident and determination of its existence
actually requires a proof. This, moreover, inaugurates a new
ontology that does not identify the signifier with the signified.
Thus the square becomes a mental object and it is derived only by
logical construction, its approximative representation not being
enough. Indeed, the presentation and the proof of the Pythagorean
Theorem (as it is given in Euclid's Elements) assumed, besides
the configuration of an orthogonal triangle, the pattern of a
square and its representations as well (Elements I, 46-48).

The
previous statement suggests that the transition from the
archetypal formulation to a rigorous geometric interpretation
demanded not only a high degree of abstraction, intuition, and
invention, but also a successful work on the foundations.
Considering the epistemological character of the particular
result, we notice that a purely perceptual category such as
verticality:

  (i) Has being objectified at an early stage, through arithmetical relations (namely, Pythagorean Triads), and finally,

 (ii) Transformed into a general mathematical form (Pythagorean theorem), obtaining a universal conceptual validity.
\par
\bigskip
\centerline{\bf 4.2  Similarity}
\bigskip
Some archetypal results may constitute the first efforts of the
determination of the shape by means of numerical relationships.
In particular, archetypal results concerning right angles evolved
into the Pythagorean Theorem. On the other hand, there is
evidence that attempts for analogous numerical expressions, not
involving in certain cases only right angles appeared in ancient
Egyptian texts, under the name of "se-qet" (Heath, 1975, pp. 126
- 128). We should note that, the geometric significance of these
two results rested upon a main feature of the shape, namely the
angle.

Furthermore, the plotting of a figure under scale usually
preserves angles. Primal mathematical activities, which
traditionally are attributed to Thales, where realized on the
ground of angle invariance and they were obviously related to the
similarity of figures, as we mean it today. The similarity
relationship classifies the figures and distinguishes the shape,
apart from the magnitude. This leads to wider and directly
"readable" classes of objects. The theory of similarity
constitutes the laying down of the criteria that explained and
established the invariance of the visual form. A necessary
premise for this process is the shift from relations concerning
two figures to internal relations that refer to one and the same
figure.

We should stress that the objectification of the notion
of similarity would not be achieved without the development of
the Eudoxian theory of proportions (or, equivalently, Dendekid's
theory of real numbers). In the exceptional case of the regular
polygons the similarity is self-evident: given n (a natural
number), any two regular n-polygons are similar. As we know
today, such a powerful situation is due to structural invariance
resulting from symmetry requirements for plane figures.

 Another
remarkable case of automatic similarity is traced back to the
writings of Plato and his intention to describe the rigidity of
the forms:

...rectilinear  surface is composed of triangles, and
all triangles are originally of two kinds.... Both of which are
made up of one right and two acute angles....

 One of them has at
either end of the base half of the divided right angle, having
equal sides, while in the other the right angle is divided into
unequal parts, having unequal sides...

 Now of the two triangles,
the isosceles has one form only; the scalene or unequal - sided
has an infinite number...

 Of the infinite forms we must again select the most beautiful...
Then let us choose two triangles... one isosceles, the other
having the square of the longer side equal three times the square
of the lesser side, (Plato, Timaeus, 53d-54b).

In both the above cases, Plato achieves rigidity through the
internal ratio of two sides. In the next, he refers to the
figures derived from the division of an equilateral triangle by
its height i.e., orthogonal triangles with angles 30 and 60
degrees. He also deals with isosceles orthogonal triangles
created when drawing the diagonal of a square [16].

 Considering two
triangles, their angle equality is a consequence of their side
proportion and vice-versa, a fact which definitely fails for
other polygons. As a matter of fact, the apprehension of the
general similarity definition, as it appears (in a complete form)
in Euclid's Elements [VI, Definitions 1 and 2], deserves a
careful historical investigation. We know that in evaluating such
an attempt, inherent epistemological problems are raised, which
in Euclid's foundation are hidden behind the formulation of the
5th postulate [17].

 Since the definition of similarity requires the
equality of ratios, its testing demands the potential infinite.
Therefore, in the case of similarity criteria we have an
inherent, non-finitistic, and simultaneously transcendental
character.  It is a remarkable fact, that if we concentrate on
the fundamental pattern of the triangle then the angle's equality
is sufficient for a decision on the similarity of the triangles!
After this, the triangle has been established as the basic
methodological tool for the proof procedure in synthetic
geometry. The transcendental character of similarity was implicit
in the manipulations using triangles.

\bigskip
\centerline{\bf 5. Diffusion}
\bigskip
Archetypal results are transformed, modified, and their variants
finally are greatly dispersed in quite different theories. The
manifold of features that allowed changes is related to the
intension of the result. Often, the context itself can be a
factor of change leading to the transformation of the original
result and its re-emergence in another theory. In the rest of our
study, for this process we utilize the term diffusion. Insofar as
archetypal results constitute aspects of objectification of the
perceptual categories, the proposed approach has obviously an
epistemological and cognitive character. On the other hand, their
reactivation in scientific programs and proposals demands a
purely mathematical and historical research. Our purpose is to
carry out parts of the above project related to the reactivation
of the Pythagorean Theorem, focusing on certain diffusions.

 In this respect, there is an early direct generalization [Euclid's Elements,
  VI.31], which replaces squares with arbitrary similar polygons drawn on
   the sides of an orthogonal triangle and proves that the same relation
    concerning their areas is also true. We have in this case the change of
     single element of the original result, in particular the squares, while
      the rest of the context is left invariant. Another modification is
       obtained after the variation of the right angle to an acute or obtuse
       angle [Elements, II.12 -13], where a well-known generalization is
       again proved. Along the same lines we can find corresponding results
       about areas and volumes of  'polytopes' in space [18].

In the Cartesian approach to Geometry, we read the Pythagorean
Theorem as a result involving magnitudes, representing directly
the sides of the orthogonal triangle: the knowledge of two of
them allows the calculation of the third, $ a=\sqrt{b^{2} +
c^{2}}$.
  This formulation
involves not directly the geometric objects, but the mediation of
an algebraic representation for them and actually constitutes in
a conceptual shift. We recognize that the obtained relation is
not given in terms of 'first reference' (Klein, 1981, p, 28).
Thus, the Pythagorean Theorem becomes now fundamental, as a
calculating tool involving lengths of line segments, in terms of
coordinates. Furthermore, this establishes the quadratic form for
the metric, a fact of paramount importance for modern mathematics.

In this modern perspective, a new insight about space conception
will lead to non-Euclidean Geometry. In such deviating contexts,
exact formulations (like Hyperbolic Trigonometry) are again
available, which approximately resemble the original Pythagorean
Theorem. At this stage of scientific evolution, an
epistemological evaluation of the original result is possible. It
is not surprising, inasmuch as we know today that this theorem is
in fact equivalent to the 5th Postulate and to the Theory of
Similarity as well [19]. Thus the incorporation and the proof of an
archetypal result into a theory, probably acts as a pivot that
also dominates the constitution of the whole theory.

 After the
development of the Infinitesimal Concepts and Calculus,
'approximation' methods were adopted by Geometry. Differential
Geometry realizes a distinction upon notions of local or global
as fields of investigation. In this setting, the metric character
inherited from the Pythagorean Theorem has survived, although the
Euclidean situation refers to the first steps of an approximating
procedure, determining the local [20]. However, variations of the
original Pythagorean result can be traced in almost every modern
mathematical theory, even if they are not directly related to the
classical views of Geometry. Along these lines, Kandisson (2002)
presented interesting implications for the extension and
utilization of the Pythagorean Theorem.

 All the above suggest a
diffusion of archetypal results in the total corpus of
Mathematics. At the same time, they constitute a confirmation of
the objectification process of our basic perceptual categories
and the reactivation of the idea of embodiment in the more
abstract mathematical fields.
\par
\bigskip
\centerline{\bf 6. Conclusions}
\bigskip
The main goal of our study is to make a contribution to the
theory of embodied mathematics and its epistemological
consequences. The main view of the authors is guided by the
interpretations of the phenomenological tradition. In the
Epistemology of Mathematics, this imposes the enlargement of our
horizons in order to include human experience and its inherent
possibility for transformations (Varelas \& al). The very
mechanism is based on experience, but could not be reduced solely
to a psychological basis. According to Vygotsky (1988, 74), "all
higher mental functions are internalized social relationships".
The biological subjects participate in historical events and
interweave into the origins of Geometry [21].

 We begin with the
perceptual categories of verticality, horizontality and
similarity (: recognition of the shapes and the angles). Humans
turned to the objective determination of these categories by
using numerical relationships and created archetypal results. In
the Greek period of mathematics, an entirely new context for
objectification was set up, based on Logic and the requirements
of the anti-visual and the anti-empirical. Thus the ideal
concepts (: prototypes) of point, line, plane, etc were invented
and proofs were achieved.

 Perceptual categories are transformed
into conceptual formation, Tall \& al (2000a). In Geometry, this
modification is expressed in terms of the notions of
perpendicular, parallel, and similar. In order to describe the
archetypal results in this new environment, new tools are needed
and the basic pattern of the triangle is created. Finally, the
Pythagorean Theorem and theorems relevant to similarity were
proved in Euclid's Elements, constituted the main axes of
Geometry, contributed to the logical construction of space.

 We
summarize the previous presentation in Table 1, specifying the
terminology and indicating the main steps of the program.
\par
\bigskip
\centerline{Table 1}
\bigskip

\centerline{\bf Notes}
\bigskip
[1] Similar questions and answers occurred in Lakoff (1987).

[2] Our civilization is geometric. The expression comes from
Freudental (Beth - Piaget, p. 222).

[3] "the phenomenal orientation of the form is determined by
directions in environment. These directions are supplied by the
pull of gravity, the visual frame of reference, or instructions",
Zuzne L. (1970), Visual Perception of Form, Academic Press, p.
301.

[4] "All effects of constancy, including that of Gestalt, are
based on the single function of extricating the essential factor
by abstracting from the inessential sensory data. The
differentiation of this function attains an amazing development
in service of shape constancy, and it needs only to be driven one
little step further to make possible an absolutely new operation
miraculously analogous to the formation of abstract, generic
concepts. Not only the small children, but also higher birds and
mammals, are able to perceive a supra-individual, generic Gestalt
in all the individual objects of the same kind. I hold that
Gestalt perception of this type is identical with that mysterious
function which is generally called "Intuition", and which is
generally called cognitive faculties of man. When the scientist,
confronted with a multitude of irregular and apparently
irreconcilable facts,  suddenly 'sees' the general regularity
ruling them all, when the explanation of the hitherto
inexplicable all 'at once' jumps out at him with the suddenness
of the revelation, the experience of this happening is
fundamentally similar to that other when the hidden Gestalt in a
puzzle-picture surprisingly starts out from the confusing
back-ground of irrelevant detail. The German expression in die a
Augen springen (to spring to the eyes) is very descriptive of
this progress", Smith M., (1966) Spatial Ability, University of
London Press, p. 216). See also, Lehar S. (1999), Gestalt
Isomorphism and the quantification of Spatial Perception, Gestalt
Theory, p. 133.

[5] This perhaps constitutes the first step towards the
mathematization of human experience, according to the paradigm of
Pythagorean theory of music that concerns the echo's shapes (the
first successful non-trivial reduction of quality to quantity
obtained, Koestler A. (1968), The Sleepwalkers, Hutchinson of
London.).

 "Consider a set of elements ABC. The child may classify them according
  to their qualitative resemblance, for example, color, size and shape.
   In order that these classificatory relationships be translated into
    numerical ones, the child has to abstract from these qualities,
     so that two elements are treated at the same time as being equivalent,
      and as being different, that is as standing in serial
       relationships.", W. Mays (in Piaget, 1972, p. 5).

[6] Kneale W. \& Kneale M. (1962, The development of Logic, Oxford,
p. 2) we have the confirmation "that the notion of demonstration
attracted attention first in connection with geometry".  A. Szabo
has the opinion that the first deductive proof started with Zeno
and later went to geometry.

[7] Looking at number as a philosophical term, we may view it as a
specialization of the Aristotelian category of quantity. But,
quantity in general is available for statements of the type 'more
or less', i.e., a primary form of the awareness of quantity that
we observe children (and early civilizations as well), as it is
described by Piaget (1969) and recently Stavy R. \& Tirosh D.
(1996), Intuitive rules in science and mathematics: the case of
"more of A - more of B" Int. J. Sci. Educ. 18 (6), 653-667.

[8] "It should fit copies of itself, it should fit into the
gravitational field, it should fit into the human hand...From
these considerations we can derive the norm for bricks: A brick
must have parallel plane sides, these pairs being orthogonal to
each other. Thus the concept of the Quader has been generated by
ideation, involving, however, the concepts of the 'plane',
'parallel', and 'orthogonal'...Thousands of years of practice
have proved this form of brick to be most expedient one for the
purpose of constructing walls", (Bender \& Schreiber (1980, p.
61).

[9] "The idea that there is some core part of meaning that is
invariant across all contexts or instances of a category offers a
useful solution to this problem in principle, but in practice,
cores for many words may be difficult or impossible to identify,
just as were defining features...For instance, that the meaning
of the word line is subtly different in each of many different
contexts (e.g., 'tanding in line', 'crossing the line',  'typing
a line of text' and that the variants are constructed at the time
of hearing/reading the word from some core meaning of the word in
the combination with the context in which it occurs", Malt (1999,
p 333).

[10] "The central pre-Socratic concept of 'essence', 'ousia'

Essences Are Substances, Essences are Forms, Essences Are Paterns
of Change ... the theory of essences fits together with the
classical theory of categories, which goes back to Aristotle. In
classical theory, a category is defined by a set of necessary and
sufficient conditions: a list of inherent properties that each
member has. A definition is a list of properties that are
necessary  and  sufficient  for something to be the kind of thing
it is....
 Euclid brought the folk theory of essences into mathematics
  in a big way. He claimed that only five postulates characterized the
  essence of plane geometry as a subject matter.

   He believed that from this
   essense all other geometric truths could be derived by deduction - by
    reason alone! From this came the idea that every subject matter in
    mathematics could be characterized in terms of an essence - a short
    list of axioms, taken as truths, from which all  other truths about the
     subject matter could be deduced.", Lakoff \& Nunez, (2000, pp. 107-109).

[11] The effort of Hjelmslev (1923) is interesting, producing a
natural geometry where we have no ideal notions of line or point
(R. S.  Tragesser (1984),  Husserl and Realism in Logic and
Mathematics, Cambridge, p. 98). In this case the results as such
that the circle and its tangent have common not a point but a
small arc.

[12] Definition is not a matter of giving some fixed set of
necessary and sufficient conditions for the application of a
concept; instead, concepts are defined by prototypes and by types
of relations of prototypes. Rather than being rigidly defined,
concepts arising from our experience are open - ended, Lakoff \&
Johnson (p. 125)

[13] The use of the terminology 'similar' for angles is common
(according to Proclos) in Thales  (Heath History of Greek
Mathematics, Ch. 4, 4.b) and is alive even in the definitions of
solid geometry (Euclid's Elements XI, definition 10).

[14] Today in Mathematics the notions of the angles and their
measurement is a deep topological result, (Dieudonne appendix II).

[15] They are special cells in the perception of angles. The
recognition of the angles is a sort of innate cognitive
apparatus, as evidence see in R. N. Haber \& M. Hershenson (1974),
The psychology of Visual Perception, Holt, Rinehart and Winston,
London, New york.  " ... cells have been found which respond to
the angles between two lines, rather than to the lines alone",
(p. 55) and about the infant's perception of angles (p. 358). See
also in Wenderoth P. \& D. White (1979), Angle - matching
illusions and perceived orientation, Perception Vol. 8, pp.
565-575.

[16] A detailed analysis of Plato's geometric ideas, in connection
with the above situations, can be found in Popper's writings in
particular see his essay "Plato and Geometry" pp. 251-270 in K.
Popper (1988), The world of Parmenides, Routlege.

[17] After the discovering of  'Non-Euclidean Geometry', the
Theory of Similarity constitutes a characterization of 'Euclidean
Geometry', since it is proved that the relations of 'similarity'
and 'congruence' are distinct if and only if geometry is
Euclidean.

[18] There is an extensive bibliography for n-dimensional
extensions of the Pythagorean Theorem. All this and recent
results can be found in the article of J. P. Quadrat \& J. B.
Lassere \& J. B. Hiriant - Urruty, Pythagora's Theorem for Areas,
in the American Mathematical Monthly, 108, 2002.

[19] The proofs
of these facts were carried out through the development of
Hyperbolic Geometry and its models. Thus these are 'metatheorems'
in Elementary Geometry. However, there are proofs which do not
refer directly to the models, see for instance, Millman and
Parker: Modern Geometry, A Metric approach with Models,
Springer-Verlang, 1981, pp. 219 -227.

[20] D. Laugwitz (1999), "B.
Riemann, Turning point in the conception of Mathematics",
Birkhaeuser, Boston.

[21] "Geometry and the science most closely
related to it have to do with space-time and the shapes, figures,
also shapes of motions, alternations of deformation, etc., that
are possible with space-time, particularly as measurable
magnitudes. It is now clear that even if we know almost nothing
about the historical surrounding world of the first geometers,
this much is certain as invariant, essential structure: that it
was a world of 'things' (including the human beings themselves as
subjects of this world); that all things necessarily had to have
a bodily character - although not all things could be mere
bodies, since the necessarily coexisting human beings are not
thinkable as mere bodies and, like even cultural objects which
belong with them structurally, are not exhausted in corporeal
being,"  E. Husserl (1999,  p 375).
\par
\bigskip
Ackowledgements: This paper is part of a research program funded
by The Univeristy of Athens, No 70/4/4921.
\par
\bigskip
\centerline{\bf References}
\bigskip
E. W. Beth - J. Piaget (1966): Mathematical Epistemology and
Psychology, D Reidel P.C., Dordrecht - Holland.

 Bender P. \&
Schreiber A. (1980) The Principle of Operative Concept Formation
in Geometry Teaching, Educ. Studies in Mathematics 11, pp. 59-90.

 Boero P. \& Bazzini L. \& Garutti R, (2001) Metaphors in
teaching and learning Mathematics: A case study concerning
inequalities, International Commission for the Study and
Improvement  of Mathematics Education CIEAEM 53 Verbania-Italy
21-27 (July).

Derrida J. (1962),  "Introduction" a Edmund
Husserl, L' Origin de la geometrie. Paris.  Eng. Trans. J. P.
Leavey. Stony Brook, NY, 1978.

Dieudonne (1972), Linear Algebra and Geometry, Kershaw, London.

Drodge N. E. \& Reid D. (2000), Embodied Cognition, Mathematical
Emotional Orientation, Mathematical thinking and learning, 2(4),
pp. 249-267, Lawrence Erlbaum Associates, Inc.

 Gray, E. M. \&
Tall, D. O. (2001), Relationships between embodied objects and
symbolic precepts: an explanatory theory of success and failure
in mathematics. PME25.

Gentner D. \& Medina J. (1998), Similarity
and the development of rules, Cognition, 65 (1998) pp. 263-297.

Edwards L (1998), Embodying Mathematics and Science: Microworlds
as Representations, Journal of Mathematical Behaviour 17 (1) pp.
53-78.

 T. L. Heath (1981), A History of Greek Mathematics, Oxford
(Rep. Dover 1975). New York.

T. L. Heath (1956), The Thirteen
Books of Euclid's Elements, Dover, New York.

 E. Husserl (1982),
Logical Investigations, Rootlege \& Kegan Paul, London. E. Husserl
(1998), The Crisis of the European Sciences and Transcendental
Phenomenology, Northestern University Press, Evanston.

Ibbotson A. \& Bryant P. (1976), The perpendicular Error and the
Vertical Effect in Cildren Drawing,  Perception, Vol. 5, pp.
319-326.

 Kadison V. R. (2002), The Pythagorean Theorem: I. The
finite case, Proc. Nat. Acad. Sciences, Vol. 99, no. 7 pp. 4178-84.

 J. Klein (1981), The World of Physics and The "Natural"
World, St. Johns Review, Vol. Autumn  pp. 22-34.

 Lakatos I.
(1997),  Mathematics, Science and Epistemology, Philosophical
papers Vol. II, (Edit. Waren J. \& Gregory . C), Cambridge Edit.,
London - New York.

 Lakoff G. \& Johnson M. (1980), Metaphors We
live by, The University of Chicago Press, Chicago.

 Lakoff G.
(1987), Women, Fire and Dangerous Things, The University Chicago
Press, Chicago.

Lakoff G. \& Nunez E. R. (2000) Where Mathematics
comes from, basic books, New York.

Mackay C. K. \&  Brazendale A.
H. \& Wilson L. F. (1972), Concepts of Horizontal and Vertical,
Developmental Psychology, Vol. 7, No 3, pp. 232-237.

 Malt B. C.
(1999), Word meaning, in A Companion to Cognitive Science, (Ed)
Bechtel W. \& Graham G. Blackwell, pp. 331-337, Massachusetts.

M. Merleau-Ponty, (2000) Phenomenology of Perception, translation
Colin Smith, Routledge and Kegan Paul, London.

 I. Mueller (1981),
Philosophy of Mathematics and deductive structure in Euclid's
Elements" Cabridge, Massachussets - London, MIT.

 Nunez E. R,
Eduards L.\& Matos J. P (1999), Embodied cognition as grounding
for situatedness and context in mathematics education,
Educational Studies of Mathematics, 39, pp. 45-65.

 Piaget J.  \&
Inhelder B. (1956) The Child's Conception of space, Routledge \&
Kegan Paul, London. Piaget J.  (1969) The Child's Conception of
Number, Routledge \& Kegan Paul, London.

 Piaget J. (1972), The principles of Genetic Epistemology, Routledge, London.

Presmeg G. N. (2002), Mathematical idea analysis: a science of
embodied mathematics, Journal for research in Mathematics
education, Vol. 33, no 1, 2002, pp. 59-63.

 Plato, Timaeus (1989),
translation Jowett B., (pp. 1151-1212) in PLATO, The Collected
Dialogues, Bollingen Series LXXI. Princeton.

Rosh, E (1978).
Principles of categorization,  In E. Rosh \& B. Loyd (Eds).
Cognition and categorization (pp. 27-48). Hillsdale, NJ: Lawrence
Erlbaum Associates Inc.

 Schiralli M. \& Sinclair (2003), A
constructive Response to "Where Mathematics comes from",
Educational Studies of Mathematics, 00: 1-13.

 Szabo A. (1978), The Beginnings of Greek Mathematics, Reidel Pub. Com,
 Dordrecht - Holland, Boston USA.

Tall D. (1995), Cognitive Growth in Elementary and Advanced
Mathematical Thinking, Conference of international Groop for the
Psychology of Learning Mathematics, Brazil, (Vol. I, pp.171-175).

Tall D. \& Gray E, Ali B. M., Growley L., DeMarois P., McGowen M.,
Pinto M., Pitta D., Yusof Y., (2000), Symbols and the Bifurcation
between Procedural and Conceptual Thinking, Canadian Journal of
Science, Mathematics and Technology Education Behaviour, 18, 4,
pp. 1-13.

Tall D. \& Thomas M. \& Davis G. \& Gray E. \& Simpson
A. (2000a), What is the Object of the Encapsulation of a process?
Journal of Mathematical Behavior 18 (2), pp. 223-241.

Varelas F.
\& Thompson E. \& Rosh E. (1999), The Embodied Mind, MIT,
Cambridge, Massachusetts.

 Vygotsky (1988). The genesis of highes mental funcions, pp. 61-80, in Cognitive Development to Adolescence, ed. by K. Richardson - S. Sheldon, Open University.East Sussex UK.

Watson, A., (2002).  Embodied action, effect, and symbol in
mathematical growth. In Anne D. Cockburn \& Elena Nardi (Eds),
Proceedings of the 26th Conference of the International Group for
the Psychology of Mathematics Education, 4, 369-376. Norwich: UK.

Watson \& Spyrou \& Tall (2003), The relationship between physical
embodiment and mathematical symbolism: The concept of vector, (in
print) in Mediterranean Journal in Mathematics Education 2.
Cyprus.

\end{document}